\documentclass[sort&compress]{revtex4-2}
\usepackage{graphicx}
\usepackage{amsmath}
\usepackage{url,hyperref,lineno,microtype,subcaption}
\hypersetup{
    colorlinks=true,
    linkcolor=blue,
    urlcolor=cyan,
    citecolor=blue
    }
\usepackage{bm}
\usepackage{esvect}

\renewcommand\phi{\varphi}
\renewcommand\[{\begin{equation}}
\renewcommand\]{\end{equation}}
\renewcommand{\t}[1]{\mathbf{#1}}
\newcommand{\gt}[1]{\bm{#1}}

\newcommand{\dd}{\mathrm{d}}

\newcommand{\abs}[1]{\left\lvert #1 \right\rvert}
\newcommand{\norm}[1]{\left\lVert #1 \right\rVert}
\newcommand{\scalar}[1]{\left\langle #1 \right\rangle}

\newcommand{\revise}[1]{#1}

\begin{document}

\title {Frustration Propagation in Tubular Foldable Mechanisms}
\author{Adam Reddy}
\author{Asma Karami}
\author{Hussein Nassar}\email{Corresponding author: nassarh@missouri.edu}
\affiliation{Department of Mechanical and Aerospace Engineering, University of Missouri, Columbia, MO 65211, USA}

\begin{abstract}
Shell mechanisms are patterned surface-like structures with compliant deformation modes that allow them to change shape drastically. Examples include many origami and kirigami tessellations as well as other periodic truss mechanisms. The deployment paths of a shell mechanism are greatly constrained by the inextensibility of the constitutive material locally, and by the compatibility requirements of surface geometry globally. With notable exceptions (e.g., Miura-ori), the deployment of a shell mechanism often couples in-plane stretching and out-of-plane bending. Here, we investigate the repercussions of this kinematic coupling in the presence of geometric confinement, specifically in tubular states. We demonstrate that the confinement in the hoop direction leads to a frustration that propagates axially as if by buckling. We fully characterize this phenomenon in terms of amplitude, wavelength, and mode shape, in the asymptotic regime where the size of the unit cell of the mechanism~$r$ is small compared to the typical radius of curvature~$\rho$. In particular, we conclude that the amplitude and wavelength of the frustration are of order $\sqrt{r/\rho}$ and that the mode shape is an elastica solution. Derivations are carried out for a particular pyramidal truss mechanism. Findings are supported by numerical solutions of the exact kinematics.
\end{abstract}


\maketitle
\section{Introduction}
\begin{figure}
    \centering
    \includegraphics[width=\linewidth]{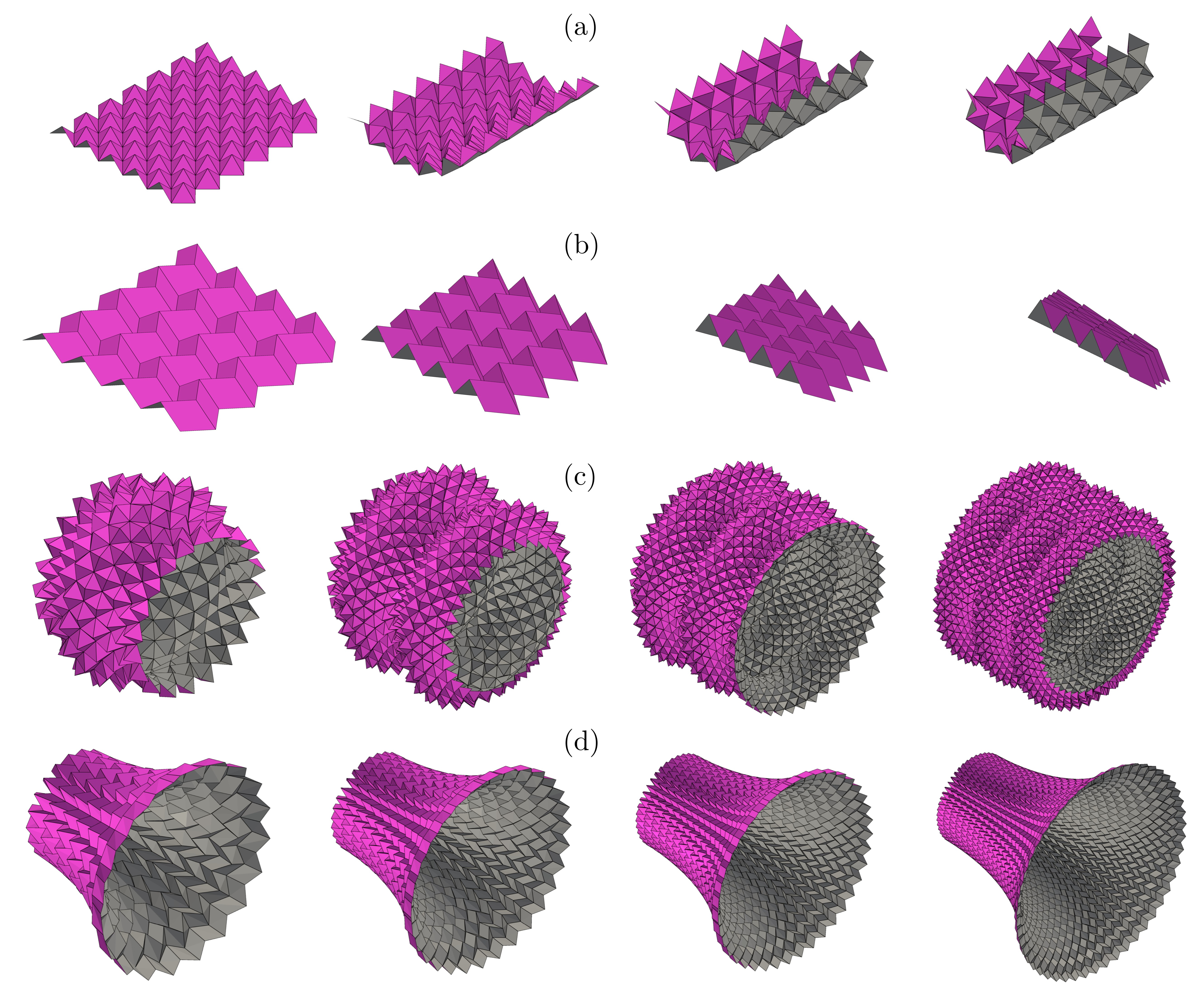}
    \caption{Coupled v.s. planar tessellations. (a, b) Uniform foldings of a coupled (a) and of a planar (b) tessellation. (c, d) Demonstration of size effects, or lack thereof: finer coupled tessellations reveal finer features (c, this paper); finer planar tessellations stabilize (d, see~\cite{Nassar2017a,Nassar2018b,Nassar2022}). \revise{Note that (c) and (d) feature non-uniformly folded states, i.e., states where folding angles vary in space.}}
    \label{fig:CoupledvsPlanar}
\end{figure}

Origami tessellations are sheets of paper folded following a repetitive crease pattern. Beyond their artistic value, they now constitute a prototyping platform for a remarkable array of foldable and deployable structures useful in engineering~\cite{Schenk2011a}. Origami tessellations are often modeled as linkages: assemblies of rigid elements representing the facets of the origami, hinged along edges representing the crease lines. From this point of view, the kinematics of origami are no different than for truss mechanisms. Thus, more generally, we will call \emph{tessellation} any linkage that is periodic (i.e., that is invariant by translation along a two-dimensional lattice of vectors) and we will call \emph{folding} any inextensional deformation that is compatible with the linkage kinematics. Suppose now that the tessellation undergoes a \emph{uniform} folding, i.e., a folding where fold angles remain periodic even if the linkage does not. Then, generically, the tessellation embraces a curved, often cylindrical, midsurface. See the examples of the waterbomb tessellation, the Yoshimura pattern and the Ron-Resch pattern~\cite{Tachi2015} as well as the pyramidal truss mechanism illustrated on Figure~\ref{fig:CoupledvsPlanar}a. Hence, it appears that the in-plane stretch, or contraction, due to the folding of a single unit cell, is coupled to an out-of-plane bending motion that is necessary to maintain geometric compatibility among consecutive unit cells. We refer to tessellations that conform to this description as \emph{coupled}. Exceptions exist however: there are tessellations that fold uniformly while maintaining a planar midsurface; we refer to them as \emph{planar}. Examples of planar tessellations include the Miura-ori (Figure~\ref{fig:CoupledvsPlanar}b), the eggbox pattern and Chebyshev nets.

Planar tessellations have received considerable attention in recent years and have been studied in the context of ``geometric mechanics'' by several authors. Schenk and Guest~\cite{Schenk2013} and Wei et al.~\cite{Wei2013} computed the Poisson's coefficient of the Miura-ori and showed that it is equal to the ratio of normal curvatures observed in its anticlastic bending. Similar results were then obtained for the eggbox pattern~\cite{Nassar2017a}, for the ``morph'' pattern~\cite{Pratapa2019}, and, more recently, for a whole class of smooth and polyhedral surfaces of translation with straight or curved creases~\cite{Nassar2022,nassar2023effective,nassar2023isometric,karami2023curvedcrease}. Nassar et al. systematically leveraged the Poisson's coefficient identity, and classical differential geometric tools, to compute the tubular states folded out of the Miura-ori~\cite{Nassar2018b}, the eggbox pattern~\cite{Nassar2017a,Nassar2017e} and the ``Mars'' tessellation~\cite{Nassar2022}. Other origami tubes were designed and investigated by Tachi and Miura~\cite{Miura2012} and by Filipov et al.~\cite{Filipov2015} with focus on axial deployment. However remarkable, the foregoing results can hardly be extended to coupled tessellations mainly because of one difficulty: size dependence. \revise{Indeed, the midsurface geometry of a coupled tessellation greatly depends on the size of the creases: finer crease patterns produce tighter midsurfaces, with divergent curvatures; see, e.g., Figure~\ref{fig:CoupledvsPlanar}c. Similarly, the waterbomb tessellation, the Ron Resch pattern and the Yoshimura pattern all admit pairs of flat-unfolded and flat-folded states, but the paths from one to the other necessarily go through curved states whose curvatures are inversely proportional to the size of the creases. By contrast, the midsurface geometry of a planar tessellation, say the Miura-ori, depends far more on the folding angles than on the size of the creases: finer crease patterns rapidly converge to a limit, non-singularly curved, surface; see, e.g., Figure~\ref{fig:CoupledvsPlanar}d.}

The purpose of the present paper is to alleviate to some extent this difficulty: we propose a differential geometric framework for the characterization of the curved midsurfaces embraced by coupled tessellations under non-uniform foldings, and provide appropriate asymptotic scalings for their size-dependent geometry. The framework is relevant for tessellations where the size~$r$ of the unit cell is small compared to the typical radius of curvature~$\rho$ of the embraced surface. \revise{Analysis of the folding of tubular states in particular reveals that confinement in the hoop direction causes bulges to develop at equal intervals in the axial direction, a phenomenon here referred to as ``frustration propagation''}. This phenomenon has recently been studied numerically by Imada and Tachi~\cite{Imada2023} as a dynamical system with an area-preserving quality that explains periodicity. Here, focus is on the size-dependence of the frustration. Specifically, we show that, generically, the amplitude and wavelength of the frustration are of order $\sqrt{r/\rho}$ and that the profile of the frustration is an elastica solution. Earlier hints of the size dependent behavior of coupled tessellations can also be inferred from the boundary layers observed in~\cite{Nassar2017a} and in~\cite{Tachi2015}.

The paper begins by introducing one particular coupled tessellation in the form of a pyramidal truss mechanism. Uniform foldings and corresponding stretch-bend kinematics are explored first. Then, non-uniform foldings are described using the differential geometry of the midsurface. The key element is an equation that determines the metric of the midsurface in function of its curvatures and of a characteristic length scale. The theoretical implications are explored for tubular states and the main results regarding frustration propagation are proven. Numerical solutions of the exact kinematics are shown to match the theoretical predictions.

\section{Theory}
\subsection{Uniform foldings}
\begin{figure}[ht!]
    \centering
    \includegraphics[width=\linewidth]{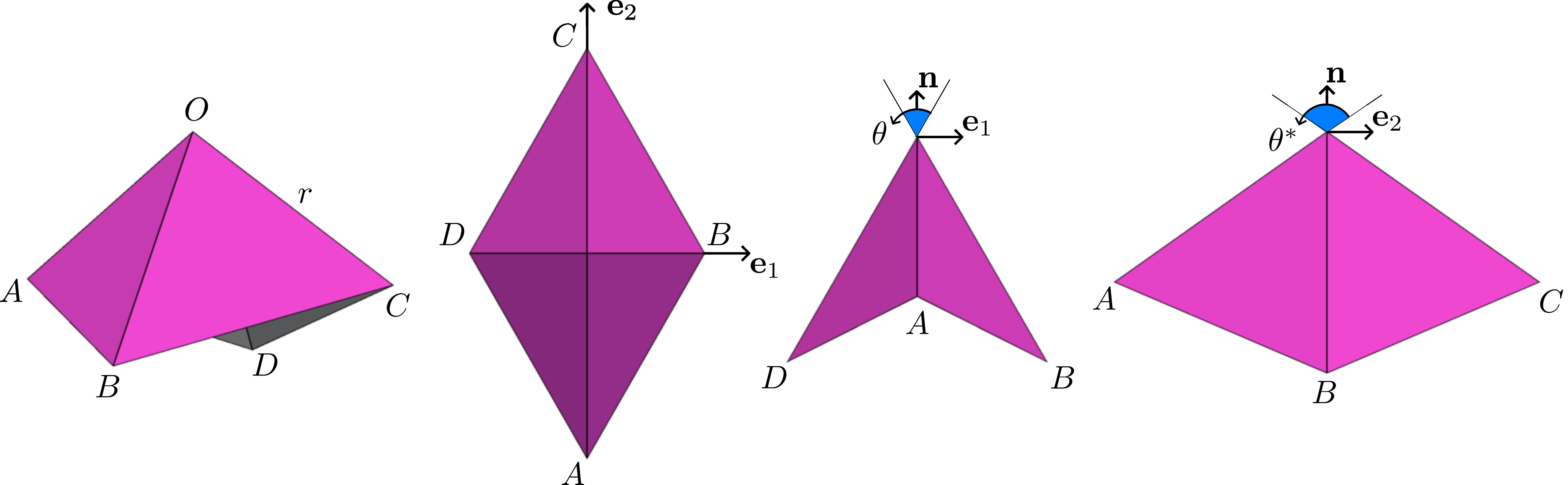}
    \caption{Annotated unit cell of the pyramidal truss mechanism: \revise{the truss is equilateral with side length~$r$.}}
    \label{fig:unitcell}
\end{figure}

\begin{figure}[t!]
    \centering
    \includegraphics[width=0.9\linewidth]{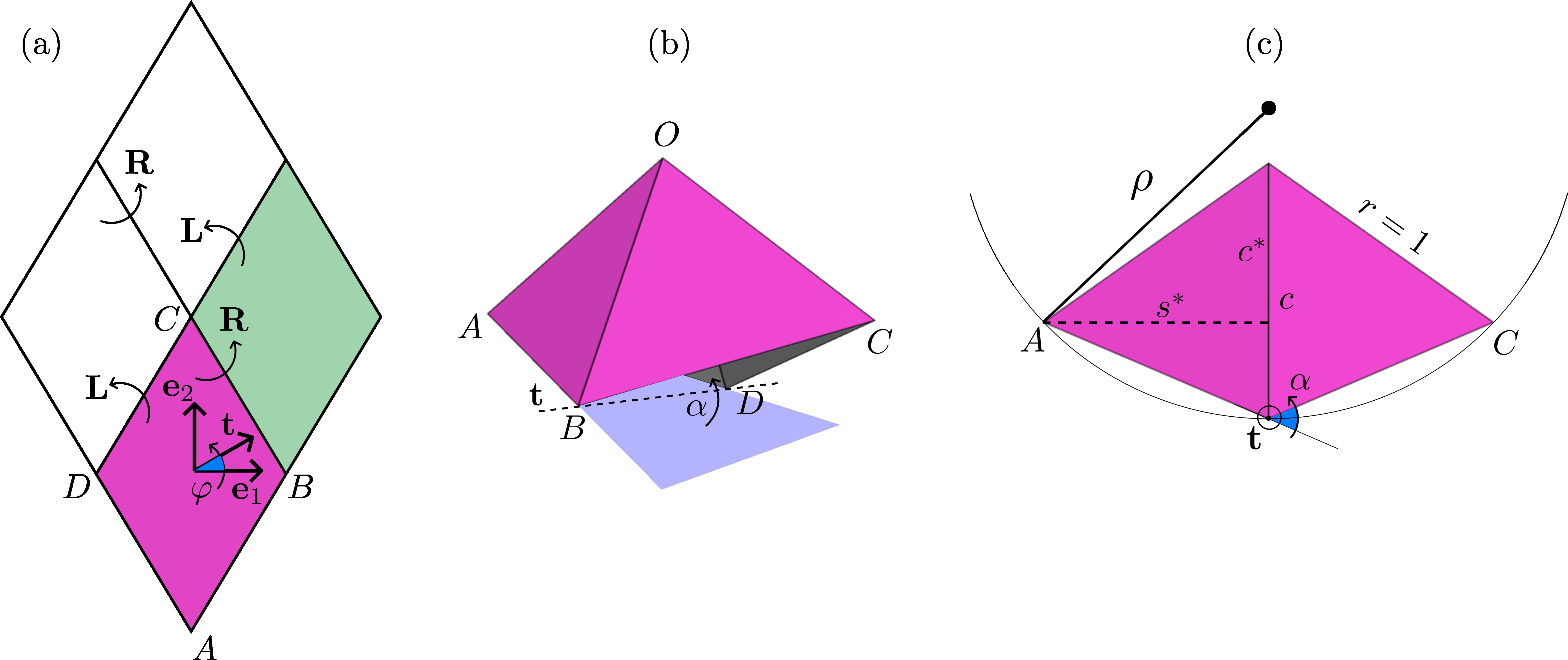}
    \caption{Analysis of uniform foldings. (a) Compatibility of rotations: \revise{arrows denote transitions between neighboring unit cells and not the axes of rotation. (b) Axis and angle of rotation in the particular case $\phi=0$. (c) Another view of (b) in the plane normal to $\t t$ showing the radius of the cylindrical state: reported angles and lengths are measured in the normal plane assuming $r=1$. For the full cylindrical state, refer to Figure~\ref{fig:CoupledvsPlanar}b.}}
    \label{fig:uniformFoldings}
\end{figure}

Consider the pyramidal truss mechanism of Figure~\ref{fig:CoupledvsPlanar}a \revise{initially introduced in~\cite{Nassar2017a}}. Its unit cell is a spherical four-bar linkage with a single Degree Of Freedom (DOF) assigned to either of the two central angles, $\theta$ and $\theta^*$; see Figure~\ref{fig:unitcell}. Let $O$ designate the apex of a pyramid with side length~$r$ and let $(A,B,C,D)$ be the vertices at its base. Let $(\t e_1,\t e_2,\t n)$ be an orthonormal basis such that $\t e_1$ and $\t e_2$ are respectively aligned with $\vv{DB}$ and $\vv{AC}$. Then, with $O$ as the origin, the remaining vertices are given by the coordinates
\[
[A] = r\begin{bmatrix}
0\\ -s^*\\ -c^*
\end{bmatrix},\quad
[B] = r\begin{bmatrix}
s\\ 0\\ -c
\end{bmatrix},\quad
[C] = r\begin{bmatrix}
0\\ s^*\\ -c^*
\end{bmatrix},\quad
[D] = r\begin{bmatrix}
-s\\ 0\\ -c
\end{bmatrix},
\]
with
\[
s \equiv \sin(\theta/2), \quad
c \equiv\cos(\theta/2), \quad
s^* \equiv \sin(\theta^*/2), \quad
c^* \equiv\cos(\theta^*/2).
\]
Last, $\theta$ and $\theta^*$ depend on one another through the constraint $AB=1$ which yields
\[
    2cc^* = 1.
\]
A uniform folding is characterized by three elements: angle $\theta$ and two rotations $\t L$ and $\t R$ that map a unit cell to two of its neighbors. Compatibility requires that $\t L$ and $\t R$ commute as shown on Figure~\ref{fig:uniformFoldings}a. Therefore, $\t L$ and $\t R$ share the same axis of rotation given by some unit vector $\t t$. Furthermore,
\[\label{eq:L&R}
    \vv{DC} = \t L\,\vv{AB}, \quad
    \vv{BC} = \t R\,\vv{AD}.
\]
Thus, compatible axes of rotation $\t t$ must be in the bisector plane of the diad $(\vv{AB},\vv{DC})$ as well as in that of $(\vv{AD},\vv{BC})$. These two planes contain $\t e_1$ and $\t e_2$ and are, in fact, the same. Therefore, the axis of rotation is given by some angle $\phi$ such that
\[
    \t t = \t e_1\cos\phi + \t e_2\sin\phi.
\]
Once $\phi$ is given, $\t R$ and $\t L$ are uniquely determined through~\eqref{eq:L&R}. By iterating $\t L$ and $\t R$, the whole folded state can be constructed out of a single unit cell. Given that $\t L$ and $\t R$ share the same axis $\t t$, and that $\t t$ is invariant under the action of $\t L$ and $\t R$, it comes that the folded state is a cylinder of axis $\t t$. In conclusion, uniformly folded states constitute a 2-DOF family of cylinders parametrized by $\theta$ and $\phi$: the former prescribes the folding angles within a unit cell whereas the latter prescribes the folding angles between neighboring cells. \revise{For other tessellations that embrace similar cylindrical states, see~\cite{Tachi2015, Feng2020}.}

\begin{figure}[ht!]
    \centering
    \includegraphics[width=0.9\linewidth]{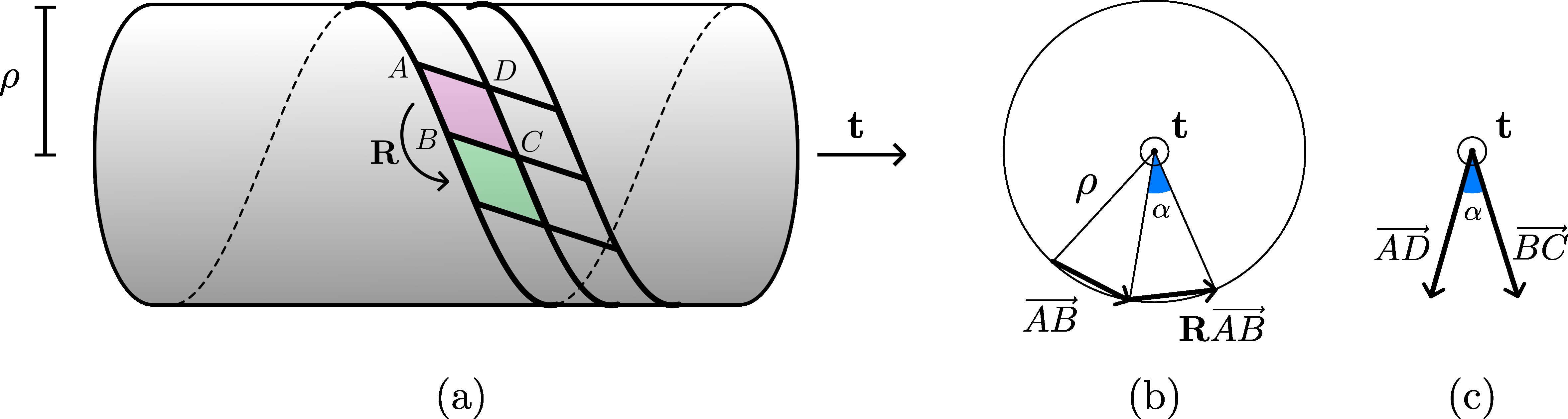}
    \caption{Geometry of a general uniformly folded state: (a) schematic view of the state; (b) axial view featuring quantities relevant to the determination of $\rho$; (c) same view featuring quantities relevant to the determination of $\alpha$.}
    \label{fig:uniformFoldingsbis}
\end{figure}

To gain further insight, let us explore the particular case $\phi=0$, i.e., where the axis of rotation $\t t$ aligns with $\t e_1$. Then, by symmetry, whatever rotation maps $\vv{AB}$ to $\vv{DC}$ also maps $\vv{AD}$ to $\vv{BC}$. In other words, $\t R=\t L$. The common angle of rotation denoted $\alpha$ can be deduced algebraically from equation~\eqref{eq:L&R} or geometrically from Figure~\ref{fig:uniformFoldings}(b, c), namely
\[
\tan(\alpha/2) = \frac{c-c^*}{s^*}.
\]
Hence, from the same Figure, it comes that the radius of the embraced cylinder, as measured from the cylinder axis to the base of a pyramid, is
\[
    \rho_{|\phi=0} = \frac{cr/2}{\abs{\sin(\alpha/2)}} = \frac r2\frac{c^2}{\abs{c-c^*}}.
\]
\revise{More generally, consider how a unit cell paves a cylinder by the repeated action of rotations as shown on Figure~\ref{fig:uniformFoldingsbis}a. By inspection of Figure~\ref{fig:uniformFoldingsbis}b, it comes that the radius of curvature is
\[
    \rho = \frac{\norm{\vv{AB}-\scalar{\vv{AB},\t t}\t t}}{2\abs{\sin(\alpha/2)}}
\]
where $\alpha$ is the angle of rotation $\t R$ and is given by
\[
\cos\alpha = \frac{\scalar{\vv{AD}-\scalar{\vv{AD},\t t}\t t,\vv{BC}-\scalar{\vv{BC},\t t}\t t}}{\norm{\vv{AD}-\scalar{\vv{AD},\t t}\t t}\norm{\vv{BC}-\scalar{\vv{BC},\t t}\t t}};
\]
see Figure~\ref{fig:uniformFoldingsbis}c.} Expanding the various dot products leads to
\[
\rho = \frac{r}{2}\frac{\sqrt{1-(s^*\sin\phi+s\cos\phi)^2}\sqrt{1-(s^*\sin\phi-s\cos\phi)^2}}{\abs{c-c^*}}.
\]

\begin{figure}[ht!]
    \centering
    \includegraphics[width=0.9\linewidth]{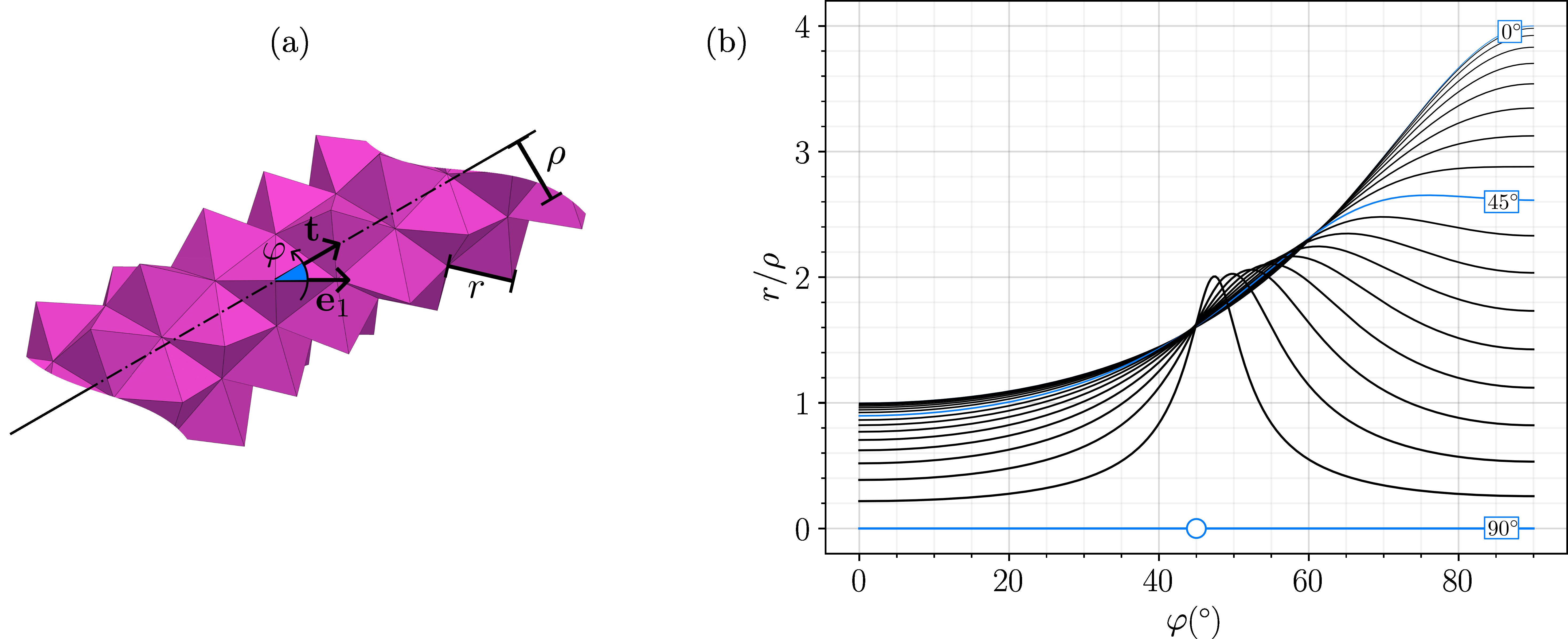}
    \caption{Variations of $\rho$ relative to $\phi$: (a) annotated generic uniform folding; (b) plot of $r/\rho$ v.s. $\phi$ for $\theta$ increasing from $0^\circ$ to $90^\circ$ by increments of $5^\circ$ (as coded by line thickness). For $\theta=90^\circ$ and $\phi=45^\circ$ (marked with a circle), curvature is indeterminate and can take any value smaller than 1.}
    \label{fig:uniformFoldingsbisbis}
\end{figure}

A generic uniform folding is illustrated on Figure~\ref{fig:uniformFoldingsbisbis}a. Profiles of the normalized curvature $r/\rho$ v.s. angle $\phi$ are depicted on Figure~\ref{fig:uniformFoldingsbisbis}b for various folding angles $\theta$. The curvature is smallest (globally or locally) at $\phi=0$ and at $\phi=\pi/2$; that is: the least curved states are cylinders of axes $\t e_1$ or $\t e_2$. Furthermore, the radius of curvature $\rho$ becomes comparable to $r$ as soon as $\theta$ departs significantly from $\pi/2$. For $\theta=\pi/2$, there is a singularity: in that case, all $\phi\neq \pm\pi/4$ lead to the same planar state with zero curvature. As for $\phi=\pm\pi/4$, the axis of the cylinder aligns with either side of the base and radius $\rho$ can take any value larger than $r$.

\subsection{Linearization}
The derivations of the above section show that any extension in direction $\t e_1$ or $\t e_2$, as measured by angle~$\theta$, is coupled to bending about some axis $\t t$ with a radius of curvature $\rho$. Furthermore, $\rho$ scales like $r$. In other words, for equal $\theta$, finer tessellations wind more tightly. \revise{The divergence of the curvature in $1/r$ towards infinity as $r$ decreases can only be avoided by limiting the magnitude of the folding angle $\theta$ to values close to~$\pi/2$.} Thus, henceforth, focus is on states such that
\[
    \theta = \pi/2 - \gamma,\quad
    \theta^* = \pi/2 + \gamma,\quad
    \gamma \ll 1,
\]
where $\gamma$ represents a small relative extension. The idea behind this restriction is that it enforces $r\ll\rho$ in a way that makes it possible to speak of a smooth midsurface for more general, non-uniformly folded, states. However, it is not clear how small $\gamma$ should or will be compared to $r/\rho$ and both are kept as small independent parameters for the time being. Then, to leading order, the radius of curvature is
\[
\rho = \frac{\abs{\cos(2\phi)}}{2\sqrt 2}\frac{r}{\abs{\gamma}}.
\]

Beyond $\rho$, it will prove convenient to compute how the midsurface bends relative to directions $\t e_1$ and $\t e_2$. Specifically, let
\[
E \equiv \frac{\scalar{\t L^T\t R\t e_1-\t e_1,\t n}}{DB},\quad
G \equiv \frac{\scalar{\t L\t R\t e_2-\t e_2,\t n}}{AC},\quad
F \equiv \frac{\scalar{\t L\t R\t e_1-\t e_1,\t n}}{AC}.
\]
Therein, $E$ is the normal change in $\t e_1$ transported to the next cell over in direction $\t e_1$ across a distance $DB$; thus, it quantifies the normal curvature of the midsurface in direction $\t e_1$; refer to Figure~\ref{fig:uniformFoldings}a for notations. Similarly, $G$ quantifies the normal curvature in direction $\t e_2$ and $F$ quantifies the torsion of the midsurface. The angles of rotations $\t L$ and $\t R$, called  $\alpha$ and $\beta$ respectively, are given by equation~\eqref{eq:L&R} and read
\[
\alpha = \frac{2\gamma}{\cos\phi - \sin\phi},\quad
\beta = \frac{2\gamma}{\cos\phi+\sin\phi},
\]
to leading order in $\gamma$. Therefore,
\[\begin{split}
E &= \frac{\alpha-\beta}{r\sqrt2}\sin\phi = \frac{2\gamma\sqrt2}{r}\frac{\sin^2\phi}{\cos(2\phi)},\\
G &= \frac{\alpha+\beta}{r\sqrt2}\cos\phi = \frac{2\gamma\sqrt2}{r}\frac{\cos^2\phi}{\cos(2\phi)},\\
F &= -\frac{\alpha+\beta}{r\sqrt2}\sin\phi = -\frac{2\gamma\sqrt2}{r}\frac{\cos\phi\sin\phi}{\cos(2\phi)}.
\end{split}
\]
It appears then that $(Er/\gamma,Gr/\gamma,Fr/\gamma)$ depend solely on $\phi$ meaning that these quantities must satisfy two algebraic constraints, namely
\[
G - E = \frac{2\gamma\sqrt2}{r},\quad
EG - F^2 = 0.
\]
These algebraic constraints define the accessible bent states as prescribed solely by the angular extension~$\gamma$, regardless of $\phi\neq\pm\pi/4$. The latter constraint, in particular, says that the Gaussian curvature vanishes, which is expected since cylinders have zero Gaussian curvature. The case $\phi=\pm\pi/4$ has been avoided above but, should it arise, the restriction $r\ll\rho$ implies $\gamma=0$ and both algebraic constraints remain valid.

Now since these states all reside in the vicinity of the planar state (i.e., $\theta=\pi/2$) and are accessible, a linear superposition should also be accessible, even if non-uniformly folded. For instance, by combining a state A bent about $\t e_1$ (i.e., $\phi_A=0$ and $\gamma_A\neq 0$) with a state B bent about $\t e_2$ (i.e., $\phi_B=\pi/2$ and $\gamma_B\neq 0$), a doubly curved state is obtained that is extended through $\gamma = \gamma_A+\gamma_B$ and whose curvatures are given by $E=E_B$, $G=G_A$ and $F=0$. Evidently, such a combined state will no longer be cylindrical or even have zero Gaussian curvature (i.e., $EG-F^2\neq 0$) but, remarkably, it will still satisfy the first constraint, namely $G-E=2\gamma\sqrt 2/r$. More generally, this constraint holds for uniformly folded states and is linear, therefore it holds for any linear combination of folded states. By contrast, the second constraint, namely $EG-F^2=0$ holds for uniformly folded states but does not carry over to more general states.

\subsection{Differential geometry of the midsurface}
We are now ready to make the scale transition in kinematic modeling from the discrete level to the continuum level. Let $(\xi_1,\xi_2)\mapsto\t x(\xi_1,\xi_2)$ parametrize the midsurface of a general, non-uniformly folded, state. The in-plane deformations of the midsurface can be quantified using the metric tensor $\gt g=g_{ij}$ with
\[
    g_{ij} = \scalar{\t x_i,\t x_j},\quad
    \t x_i \equiv \frac{\partial \t x}{\partial \xi_i}.
\]
Here, we identify
\[
    \t x_1 \equiv \frac{\vv{DB}_{|\gamma}}{DB_{|\gamma=0}} = (1-\gamma/2)\t e_1,\quad
    \t x_2 \equiv \frac{\vv{AC}_{|\gamma}}{AC_{|\gamma=0}} = (1+\gamma/2)\t e_2,
\]
which implies, to first order in $\gamma$, that
\[\label{eq:ong}
    [\t g] = \begin{bmatrix}
        1-\gamma & 0 \\ 0 & 1+\gamma
    \end{bmatrix}.
\]
Above, and in the following, it is understood that $\gamma\equiv\gamma(\xi_1,\xi_2)$ can depend on space coordinates so as to allow for the continuum description to encompass non-uniformly folded states. The only restriction in that regard is that variations in $\gamma$ should occur over length scales larger than $r$.

As for out-of-plane deformations, they are captured by the second fundamental form $\t b=b_{ij}$ where
\[
    b_{ij} = \scalar{\t x_{ij},\t n},\quad \t x_{ij} \equiv \frac{\partial^2\t x}{\partial\xi_i\partial\xi_j},\quad \t n \equiv \frac{\t x_1\wedge\t x_2}{\norm{\t x_1\wedge\t x_2}}.
\]
Here too, we identify $b_{11}\equiv E$, $b_{22}\equiv G$ and $b_{12}\equiv F$ so that the constraint
\[\label{eq:onb}
    b_{22}-b_{11} = \frac{2\gamma\sqrt 2}{r}
\]
is systematically enforced point-wise, i.e., at all positions $(\xi_1,\xi_2)$. Again, the other constraint, namely $\det\t b = 0$, only holds for uniformly folded states (i.e., $\gamma=\text{cste}$) but not in general.

Solving for the midsurface then amounts to finding surfaces whose metric and second form, $\t g$ and $\t b$, satisfy equations~\eqref{eq:ong} and~\eqref{eq:onb}. Alternatively, it is possible to combine both equations into one statement: \emph{the tessellation embraces midsurfaces whose metric depend on the second form through}
\[\label{eq:central}
    [\t g] = \begin{bmatrix}
        1 & 0 \\ 0 & 1
    \end{bmatrix}
    + r\frac{b_{22}-b_{11}}{2\sqrt 2}\begin{bmatrix}
        -1 & 0 \\ 0 & 1
    \end{bmatrix}.
\]

\revise{It is worthwhile to highlight that for several planar tessellations, including the Miura-ori, the eggbox pattern as well as other patterns with no stretch-to-bend coupling, the coefficients of the second fundamental form $\t b$ still satisfy a linear metric-dependent constraint similar to~\eqref{eq:onb} albeit one that is homogeneous, e.g., of the form $p(\gamma)b_{22}-q(\gamma)b_{11}=0$; see, e.g., \cite{Nassar2017a, Nassar2018b, Nassar2022}. Here, by contrast, the foregoing constraint exhibits a ``source'' term in $\gamma/r$ that is both metric- and size-dependent.}
\section{Results}
\subsection{Geometry of tubular states}
We define a tubular state as one that embraces a midsurface of revolution. Specifically, let
\[
    [\t x(\xi_1,\xi_2)] = \begin{bmatrix}
        \rho(\xi_2)\cos(q\xi_1)\\
        \rho(\xi_2)\sin(q\xi_1)\\
        z(\xi_2)
    \end{bmatrix},
    \quad
    0\leq\xi_1<W
\]
where $\rho$ and $z$ are to-be-determined functions, $q\equiv 2\pi/W$ and $W$ is the width of the tessellation in the flat reference state; see Figure~\ref{fig:revolution}a. Then, computing $\t x_1$, $\t x_2$ and so on, it comes that
\[\label{eq:gamma}
    g_{11} = q^2\rho^2,\quad g_{22} = \rho'^2+z'^2,\quad
    \gamma = \frac{r}{2\sqrt2}\frac{\rho''z'-\rho'z''+q^2\rho z'}{\sqrt{\rho'^2+z'^2}},
\]
where a prime denotes $\dd/\dd \xi_2$. The key equation~\eqref{eq:central} then yields two Ordinary Differential Equations (ODEs), namely
\[\label{eq:ODE}
    q^2\rho^2 = 1 - \gamma,\quad
    \rho'^2+z'^2 = 1 + \gamma.
\]
In the following, we solve these equations in the asymptotic regime that is consistent with the adopted continuum description, i.e., in the limit where $\gamma\to 0$ and $r/\rho\to 0$.

\subsection{Asymptotics of size dependence}
\begin{figure}
    \centering
    \includegraphics[width=\linewidth]{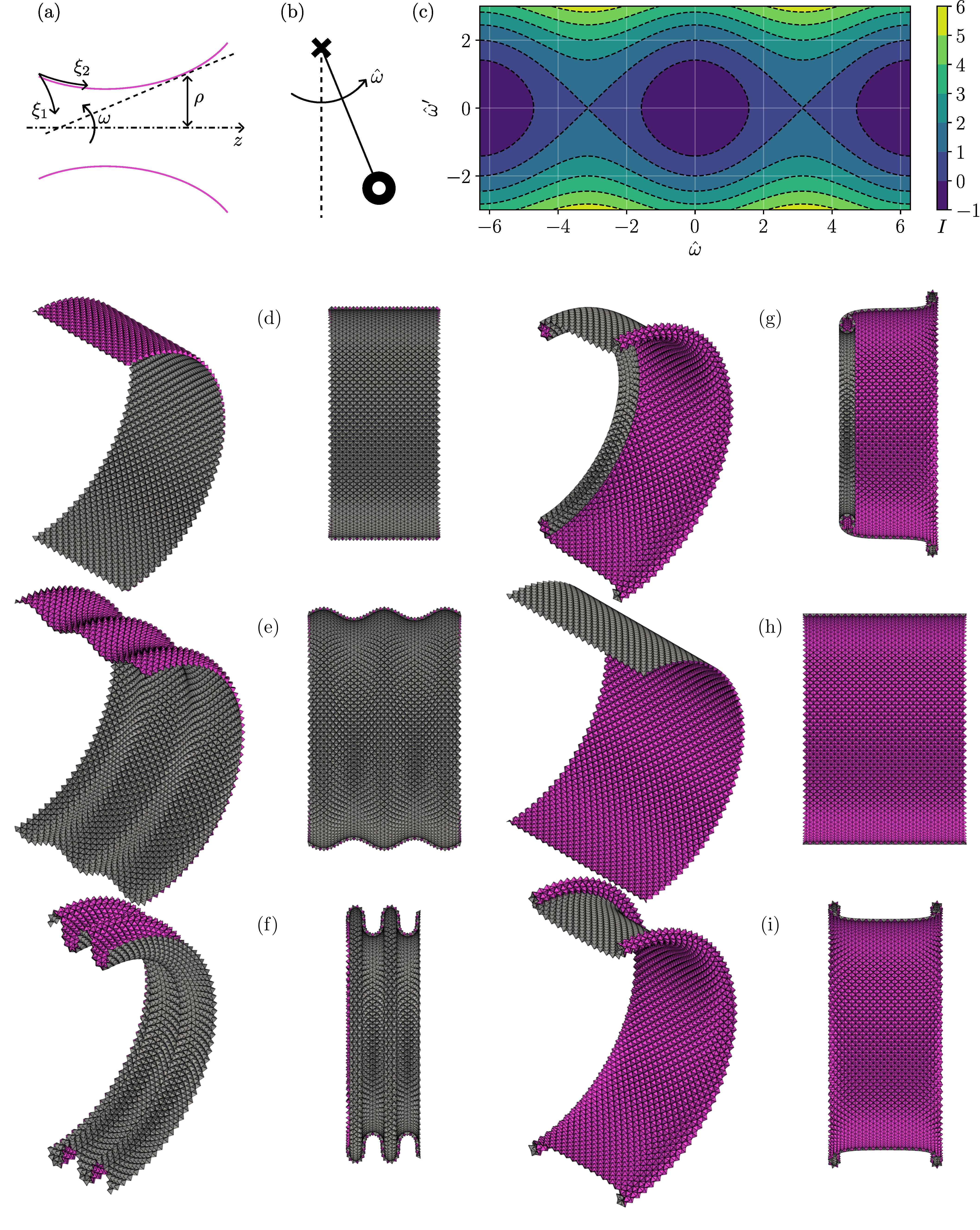}
    \caption{Tubular states. (a) Notations. (b) Simple pendulum analogy. (c) Phase diagram: iso-contours of the analogue to mechanical energy $I$. (d-i) Two views of various tubular states for increasing $I$: \revise{(d, h) are uniformly folded states; the other states are non-uniformly folded.}}
    \label{fig:revolution}
\end{figure}

We deal with equations~\eqref{eq:ODE} as algebraic equations to begin with. Then, up to an error smaller than $\gamma$,
\[
    q\rho = 1 - \frac{\gamma}{2},
\]
and there exists a function~$\omega\equiv\omega(\xi_2)$ such that
\[
\rho' = \sin\omega,\quad z' = \cos\omega,
\]
since, now to leading order in $\gamma$, we have $\rho'^2+z'^2=1$. Substituting back into expression~\eqref{eq:gamma} for $\gamma$ then yields
\[
\gamma = \frac{r}{2\sqrt2}\left(\omega'+q\cos\omega\right),
\]
where terms in $\gamma^2$ and in $\gamma qr = O(\gamma r/\rho)$ have been neglected, both being small relative to $\gamma$. Taking the derivative $\dd/\dd\xi_2$ and replacing $\gamma'$ with its leading-order expression in terms of $q\rho'$ leads to
\[
    -2q\sin\omega = \frac{r}{2\sqrt2}\left(\omega''-q\omega'\sin\omega\right).
\]
But $\omega'$ is of order $ b_{22}\ll1/r$ meaning that, to leading order, $\omega$ is solution to the second-order ODE
\[
    \omega'' + \frac{4q\sqrt2}{r}\sin\omega = 0.
\]
Therefore, it appears that $\omega$ oscillates at a length scale
\[
    L \equiv \sqrt{\frac{r}{4q\sqrt2}} = O\left(\sqrt{rW}\right).
\]
Accordingly, it is appropriate to re-scale the $\xi_2$ coordinate and introduce a function $\hat\omega$ such that
\[
    \omega(\xi_2) \equiv \hat\omega(\xi_2/L).
\]
Then, $\hat\omega$ is solution to the elastica ODE, also known as the simple pendulum ODE,
\[
    \hat\omega'' + \sin\hat\omega = 0.
\]
Note that further linearization of the elastica ODE is not warranted in general since there is no reason for $\omega$ to be small unless further assumptions regarding the smallness of the initial conditions $\hat\omega_o$ and $\hat\omega'_o$ are made. In any case, once a solution to the elastica ODE is chosen by setting $\hat\omega_o$ and $\hat\omega_o'$ at, say, $\xi_2=0$, the profile of the state can be determined by integration as in
\[
    \rho = \rho_o + \int_0^{\xi_2}\sin\hat\omega(s/L)\dd s,\quad
    z = z_o + \int_0^{\xi_2}\cos\hat\omega(s/L)\dd s,
\]
or, by re-scaling, 
\[\label{eq:rhoz}
    \rho = \rho_o + L\int_0^{\xi_2/L}\sin\hat\omega(s)\dd s,\quad
    z = z_o + L\int_0^{\xi_2/L}\cos\hat\omega(s)\dd s.
\]
By analogy with the simple pendulum (Figure~\ref{fig:revolution}b), we know that $\sin\omega$, a quantity analogous to the horizontal deflection of the pendulum, oscillates periodically with zero average. Hence, $\rho$ oscillates periodically with an amplitude of order $L=O(\sqrt{rW})$. By contrast, $\cos\omega$, a quantity analogous to the vertical deflection of the pendulum, oscillates periodically but does not necessarily have zero average. Hence, generically, $z$ grows linearly. \emph{In conclusion, tubular states exhibit a periodic frustration whose wavelength and amplitude are of order $\sqrt{rW}$ and whose profile is an elastica solution.}

Note that the initial value $z_o$ is arbitrary and inconsequential because it amounts to a rigid body translation along the axis of revolution. As for $\rho_o$, it can be obtained by reconsidering the expression of $\gamma$, namely
\[\label{eq:gammahat}
    \gamma = \frac{r\hat\omega'}{2L\sqrt2} = O(\sqrt{r/W}),
\]
where terms of order $qr = O(r/W)$ have been neglected in favor of terms of order $O(\sqrt{r/W})$. This provides an alternative expression for $\rho$, that is 
\[
    q\rho = 1 - \frac{r\hat\omega'}{2L\sqrt2},
\]
and with it an expression for $\rho_o$ in function of $\hat\omega'_o$. Note that, for the sake of consistency, the initial condition $\hat\omega'_o$ must be chosen of order $O(1)$, i.e., small relative to $L/r=O(\sqrt{W/r})$.

Last, it is worthwhile to shed some light on the normal curvatures $b_{11}$ and $b_{22}$. On one hand, $\omega=\arctan(\rho'/z')$ is the incidence angle of the tangent plane relative to the axis of revolution. Its derivative can directly be interpreted as the normal curvature in the axial direction, specifically
\[
    b_{22} = \frac{1}{L}\hat\omega'.
\]
This normal curvature blows up in the limit $r\to 0$ and is of order~$1/L=O(1/\sqrt{rW})$. On the other hand, the normal curvature in the hoop direction is
\[
    b_{11} = -q\cos\hat\omega,
\]
and remains finite, being of order $1/W$.

\subsection{Phase diagram}
The elastica ODE admits an invariant analogous to the mechanical energy of the simple pendulum
\[
 I \equiv \frac{1}{2}\hat\omega'^2 - \cos\hat\omega = \frac{1}{2}\hat\omega_o'^2 - \cos\hat\omega_o.
\]
The iso-contours of $I$ define the solution orbits in $(\hat\omega,\hat\omega')$-space, each orbit being in correspondence with a tubular state; see Figure~\ref{fig:revolution}c. Two orbits in particular, namely $I=\pm1$, degenerate into two fixed points corresponding to $(\omega=0,\gamma=0)$ and $(\omega=\pi,\gamma=0)$; these are two cylindrical states with the pyramids pointing outward and inward, respectively (Figure~\ref{fig:revolution}d,~h). Since $\gamma=0$ in both cases, it is insightful to seek the next order in its asymptotic expansion. Going back to equation~\eqref{eq:gamma} and letting $\rho$ be a constant lead to
\[
    \gamma = \pm\frac{\pi}{\sqrt{2}}\frac{r}{W} = O(r/W),
\]
where the sign is positive for $\omega=0$ and negative for $\omega=\pi$. In other words, achieving a cylindrical state requires closing the pyramids slightly (as seen in the axial direction) when the pyramids are pointing outward and opening them slightly when they are pointing inward.

Interestingly, these two fixed-points have different stability properties. The cylindrical state with the pyramids pointing out is stable: deviations lead to oscillations that grow throughout the profile at once (Figure~\ref{fig:revolution}d-f). By contrast, the cylindrical state with the pyramids pointing in is unstable: deviations lead to oscillations growing at localized intervals that, should self-penetration be precluded, cannot exist except at the edges of the tubular state (Figure~\ref{fig:revolution}g-i). These edge states correspond to a pendulum that remains near the top (unstable) equilibrium position most of the time and that, however swiftly, swings by the bottom (stable) equilibrium position.

\subsection{Numerical solutions and error analysis}
\begin{figure}
    \centering
    \includegraphics{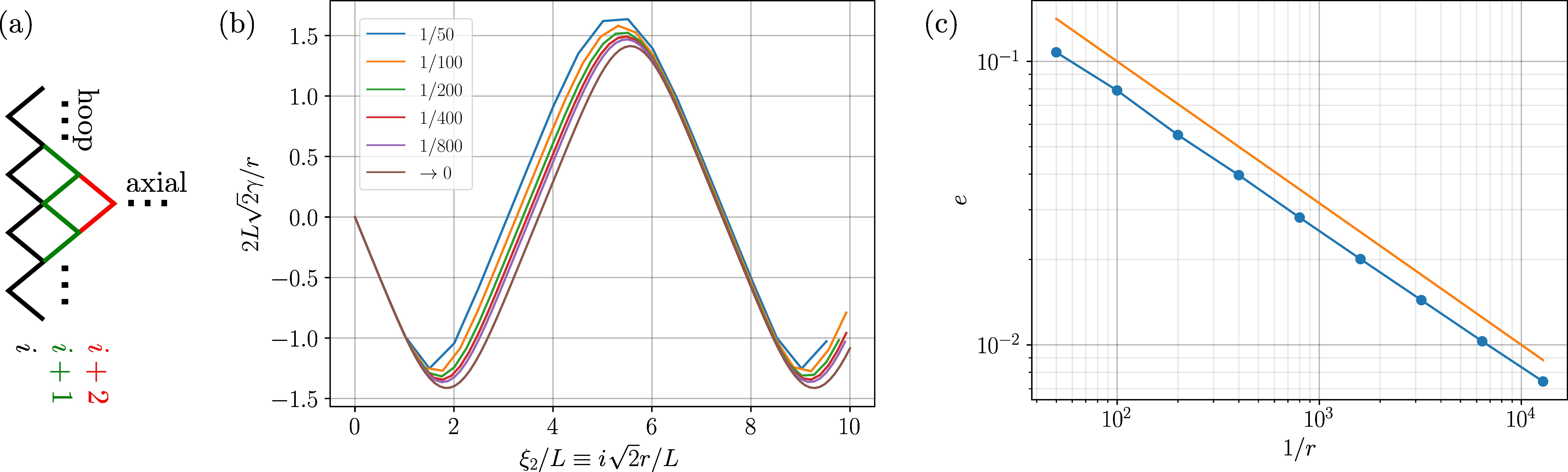}
    \caption{Convergence analysis. (a) \revise{Schematics of the discrete solution algorithm: knowledge of zigzag $i$ suffices to uniquely determine zigzag $i+1$.} (b) Normalized folding angle $\gamma$ v.s. normalized discrete curvilinear coordinate as measured by a nodal index $i$: curves are obtained for decreasing~$r$ (see legend); the limit $r\to 0$ is the solution to the elastica ODE. (c) Convergence speed: Data points show the mean quadratic error $e$ as obtained from (b); the seamless line shows the trend $e=r^{1/2}$. Analysis is carried for $\hat\omega_o=\pi/2$, $\hat\omega'_o=0$. Other states show similar trends.}
    \label{fig:cv}
\end{figure}

The tessellations produced throughout the paper and in Figure~\ref{fig:revolution} in particular are rendered by solving the exact discrete kinematics. \revise{The solution algorithm proceeds iteratively by constructing the nodes of the folded state one ``zigzag'' at a time as illustrated on Figure~\ref{fig:cv}a. Mainly, given the nodes that belong to zigzag number~$i$, the nodes of zigzag number $i+1$ can be constructed by completing the bases of the pyramids that are enclosed between zigzags $i$ and $i+1$, one pyramid at a time. This is possible because each enclosed pyramid has a unique degree of freedom.} Further detail can be found in~\cite{Nassar2017a, Nassar2018b}. The algorithm is implemented in Python and its code is available online (see Supplemental Data). Here, rotation symmetry is leveraged and the algorithm can be initialized simply by giving the three parameters $(r/W,\theta_o,\omega_o)$, where $\theta_o=\pi/2-\gamma_o$. \revise{Profiles of $\gamma$ are computed for small but finite $r$ by solving the discrete kinematics and are shown on Figure~\ref{fig:cv}b. The profile of $\gamma$ obtained by solving the elastica ODE is shown as well. The plots demonstrate that the discrete solution, hereafter denoted $\gamma_\text{discrete}$, converges to the continuous one, hereafter denoted $\gamma_\text{ODE}$. The convergence error
\[
e \equiv \frac{1}{T}\sqrt{\int_0^T\abs{\frac{2L\sqrt{2}}{r}\gamma_\text{ODE}-\frac{2L\sqrt{2}}{r}\gamma_\text{discrete}}^2}
\]
is further computed over a normalized $r$-independent range $0\leq \xi_2/L\leq T$ and is plotted on Figure~\ref{fig:cv}c. Note that the $\gamma$ profiles have been normalized as well as anticipated by equation~\eqref{eq:gammahat}. The logarithmic scale shows that the error~$e$ decays with the side length~$r$ like~$r^{1/2}$.}

\section{Conclusion}
The paper proposes a differential geometric framework that allows to deal with the size-dependent kinematics of coupled tessellations, i.e., of periodic foldable structures where in-plane folding is coupled to out-of-plane bending. The framework succeeds in predicting size effects both qualitatively and quantitatively for tubular states where hoop confinement produces a frustration that propagates axially. The convergence of the framework is somewhat slow with the error being of order $r^{1/2}$ where $r$ is a small size parameter. \revise{Thus, improvements that take higher order corrections into account are desirable as they can improve convergence speed. In principle, such corrections will intervene at two places: in the linearization of the discrete kinematics and in the asymptotics of the governing non-linear ODE.}

\revise{The main conclusion of the analysis, namely that axisymmetric frustrations are elastica-shaped periodic with wavelength and amplitude of order $r^{1/2}$, is not specific to the current pyramidal truss mechanism and should be typical of coupled tessellations that satisfy two properties: ($i$) the coupling is linear, meaning that $\t g=\t g(\t b)$ is linear (affine strictly speaking); ($ii$) the tessellation is rectangular so that the coordinate lines aligned with the axes of symmetry deform without shearing. Under these conditions, equation~\eqref{eq:ong} is generic and is sufficient for our conclusion to hold. That being said, extensions to other potentially non-linearly-coupled or oblique tessellations are of interest as well and are yet to be investigated.}

\section*{Conflict of Interest Statement}
The authors declare that the research was conducted in the absence of any commercial or financial relationships that could be construed as a potential conflict of interest.

\section*{Author Contributions}
A.R. and A.K. studied discrete and continuum kinematics. A.R. and H.N. carried out numerical simulations and wrote the paper. H.N. conceived the research.

\section*{Funding}
Work supported by the NSF under CAREER award No. CMMI-2045881.

\section*{Acknowledgments}
H.N. thanks Arthur Lebée (École des Ponts) for insightful exchange.

\section*{Supplemental Data}\label{sec:SD}
The code that produced the figures is available at \url{https://github.com/nassarh/pyramids}.

\end{document}